\documentclass[12pt]{amsart}
\usepackage[cp1251]{inputenc}
\usepackage[english]{babel}
\usepackage{amsmath}
\usepackage{amssymb}
\usepackage{amsfonts}

\begin{document}

\title[Operator symbols]
      {Operator symbols}
\author{Vladimir~B. Vasilyev}
\address{Chair of Differential Equations \\
 Belgorod State National Research University\\
         Studencheskaya 14/1\\
        308007 Belgorod, Russia}

\email{vbv57@inbox.ru}

\begin{abstract}
  We consider special elliptic operators in functional spaces on manifolds with a boundary which has some singular points. Such an operator can be represented by a sum of operators, and for a Fredholm property of an initial operator one needs a Fredholm property for an each operator from this sum.
\end{abstract}

\maketitle

\section{Introduction}
This paper is devoted to describing structure of a special class of linear bounded operators on a manifold with non-smooth boundary. Our description is based on Simonenko's theory of envelopes \cite{S} and explains why we obtain distinct theories for pseudo-differential equations and boundary value problems and distinct index theorems for such operators.

\subsection{Operators of a local type}
In this section we will give some preliminary ideas and definitions from \cite{S}.

Let $B_1,B_2$ be  Banach spaces consisting of functions defined on compact $m$-dimensional manifold $M$, $A: B_1\rightarrow B_2$ be a linear bounded operator, $W\subset M$, and $P_W$ be a projector on $W$ i.e.
$$
(P_Wu)(x)=\begin{cases}
u(x),&\text{if $x\in W$;}\\
0,&\text{if $x\notin\overline{W}$.}
\end{cases}
$$
{\bf  Definition 1.} {\it
An operator $A$ is called an operator of a local type if the operator
\[
P_UAP_V
\]
is a compact operator for arbitrary non-intersecting compact sets $U,V\subset M$.
}

\subsection{Simple examples} These are two simplest examples for the illustration.

{\bf Example 1.} {\it
If $A$ is a differential operator of the type
\[
(Au)(x)=\sum\limits_{|k|=0}^{n}a_k(x)D^ku(x),~~~D^ku=\frac{\partial^ku}{\partial x_1^{k_1}\cdots\partial x_m^{k_m}},
\]
then $A$ is an operator of a local type.
}

{\bf Example 2.} {\it
If $A$ is a Calderon--Zygmund operator with variable kernel $K(x,y)\in C^1(\mathbb R^m\times(\mathbb R^m\setminus\{0\})$ of the following type
\[
(Au)(x)=v.p.\int\limits_{\mathbb R^m}K(x,x-y)u(y)dy.
\]
}

{\bf Everywhere below we say ``an operator'' instead of ``an operator of a local type''.}

\subsection{Functional spaces on a manifold}

\subsubsection{Spaces $H^s(\mathbb R^m), L_p(\mathbb R^m), C^{\alpha}(\mathbb R^m)$}
It is possible working with distinct functional spaces \cite{MPr,T}.

{\bf Definition 2.}  \cite{E} {\it
The space $H^s(\mathbb R^m), s\in\mathbb R,$ is a Hilbert space of functions with the finite norm
\[
||u||_s=\left(\int\limits_{\mathbb R^m}|\widetilde u(\xi)|^2(1+|\xi|)^{2s}d\xi\right)^{1/2},
\]
where the sign $\sim$ over a function means its Fourier transform.
}

{\bf Definition 3.} \cite{MPr} {\it
The space $L_p(\mathbb R^m), 1<p<+\infty,$ is a Banach space of measurable functions with the finite norm
\[
||u||_p=\left(\int\limits_{\mathbb R^m}|u(x)|^pdx\right)^{1/p}.
\]
}

{\bf Definition 4.} \cite{MPr} {\it
The space $C^{\alpha}(\mathbb R^m), 0<\alpha\leq1,$ is a space of continuous on $\mathbb R^m$ functions $u$ satisfying the H\"older condition
\[
|u(x)-u(y)|\leq c|x-y|^{\alpha},~~~\forall x,y\in\mathbb R^m,
\]
with the finite norm
\[
||u||_{\alpha}=\inf\{c\},
\]
where infimum is taken over all constants $c$ from above inequality.
}

\subsubsection{Partition of unity and spaces $H^s(M), L_p(M), C^{\alpha}(M)$}

If $M$ is a compact manifold then there is {\it a partition of unity} \cite{M}. It means the following. For every finite open covering $\{U_j\}_{j=1}^k$ of the manifold $M$ there exists a system of functions $\{\varphi_j(x)\}_{j=1}^k, \varphi_j(x)\in C^{\infty}(M),$ such that
\begin{itemize}
\item  $0\leq\varphi_j(x)\leq 1$,
\item $supp~\varphi_j\subset U_j$,
\item $\sum\limits_{j=1}^k\varphi_j(x)=1$.
\end{itemize}

So we have
\[
f(x)=\sum\limits_{j=1}^k\varphi_j(x)f(x)
\]
for arbitrary function $f$ defined on $M$.

Since an every set $U_j$ is diffeomorphic to an open set $D_j\subset\mathbb R^m$ we have corresponding diffeomorphisms $\omega_j: U_j\rightarrow D_j$. Further for a function $f$ defined on $M$ we compose mappings $f_j=f\cdot\varphi_j$ and as far as $supp~f_j\subset U_j$ we put $\hat f_j=f_j\circ\omega_j^{-1}$ so that $\hat f_j: D_j\rightarrow\mathbb R $ is a function defined in a domain of $m$-dimensional space $\mathbb R^m$. WE can consider the following spaces \cite{E,MPr,T}.

{\bf Definition 5.} {\it
A function $f\in H^s(M)$ if the following norm
\[
||f||_{H^s(M)}=\sum\limits_{j=1}^k||\hat f_j||_s
\]
 is finite.

A function $f\in L_p(M)$ if the following norm
\[
||f||_{L_p(M)}=\sum\limits_{j=1}^k||\hat f_j||_p
\]
 is finite.

A function $f\in C^{\alpha}(M)$ if the following norm
\[
||f||_{C^{\alpha}(M)}=\sum\limits_{j=1}^k||\hat f_j||_{\alpha}
\]
 is finite.
}

\section{Operators on a compact manifold}

On the manifold $M$ we fix a finite open covering and a partitions of unity corresponding to this covering $\{U_j, f_j\}_{j=1}^{n}$ and
choose smooth functions
  $\{ g_j\}_{j=1}^{n}$ so that $supp~g_j\subset V_j,$ $\overline{U_j}\subset V_j$, and $g_j(x)\equiv 1$ for $x\in supp~ f_j, supp~f_j\cap(1-g_j)=\emptyset$.

{\bf Proposition 1.}
{\it The  operator $A$ on the manifold $M$ can be represented in the form
\[
A=\sum\limits_{j=1}^{n}f_j\cdot A\cdot g_j +T,
\]
where $T: B_1\rightarrow B_2$ is a compact operator.}

\begin{proof}
The proof is a very simple. Since
\[
\sum\limits_{j=1}^nf_j(x)\equiv 1,~~~\forall x\in M,
\]
then we have
\[
A=\sum\limits_{j=1}^{n}f_j\cdot A=\sum\limits_{j=1}^{n}f_j\cdot A\cdot g_j+\sum\limits_{j=1}^{n}f_j\cdot A\cdot(1-g_j),
\]
and the proof is completed.
\end{proof}

{\bf Remark 1.} {\it
It is obviously such operator is defined uniquely up to a compact operators which do not influence on an index.
}

By definition for an arbitrary operator $A: B_1\rightarrow B_2$
\[
|||A|||\equiv\inf||A+T||,
\]
where {\it infimum} is taken over all compact operators $T: B_1\rightarrow B_2$.

Let $B'_1,B'_2$ be  Banach spaces consisting of functions defined on $\mathbb R^m$, $\widetilde A: B'_1\rightarrow B'_2$ be a linear bounded operator.

Since $M$ is a compact manifold, then for every point $x\in M$ there exists a neighborhood $U\ni x$ and diffeomorphism $\omega: U\rightarrow D\subset\mathbb R^m, \omega(x)\equiv y$. We denote by $S_{\omega}$ the following operator acting from $B_k$ to $B_k', k=1,2$. For every function $u\in B_k$ vanishing out of $U$
\[
(S_{\omega}u)(y)=u(\omega^{-1}(y)),~~~y\in D,~~~(S_{\omega}u)(y)=0,~~~y\notin D.
\]

{\bf Definition 6.} {\it
A local representative of the operator $A: B_1\rightarrow B_2$ at the point $x\in M$  is called the operator $\widetilde A: B_1'\rightarrow B_2'$ such that $\forall\varepsilon >0$ there exists the neighborhood $U_j$ of the point $x\in  U_j\subset M$ with the property
\[
|||g_jAf_j-S_{\omega_j^{-1}}\hat g_j\widetilde A\hat f_jS_{\omega_j}|||<\varepsilon.
\]
}

\section{Algebra of symbols}

{\bf Definition 7.} {\it
Symbol of an operator $A$ is called the family of its local representatives $\{A_x\}$ at each point $x\in\overline{M}$.
}

One can show like \cite{S} this definition of an operator symbol conserves all properties of a symbolic calculus. Namely, up to compact summands
\begin{itemize}
\item product and sum  of two operators corresponds to product and sum of their local representatives;
\item adjoint operator corresponds to its adjoint local representative;
\item a Fredholm property of an operator corresponds to a Fredholm property of its local representative.
\end{itemize}

\section{Operators with symbols. Examples.}

It seems not every operator has a symbol, and we give some examples for operators with symbols.

{\bf Example 3.} {\it
Let $A$ be the differential operator from example 1, and functions $a_k(x)$ be continuous functions on $\mathbb R^m$. Then its symbol is an operator family consisting of multiplication operators on the function
\[
\sum\limits_{|k|=0}^{n}a_k(x)\xi^k,
\]
where $\xi^k=\xi_1^{k_1}\cdots\xi_m^{k_m} $.
}

{\bf Example 4.} {\it
Let $A$ be the Calderon--Zygmund operator from Example 2 and $\sigma(x,\xi)$ be its symbol in sense of \cite{MPr}, then its symbol is an operator family consisting of multiplication operators on the function $\sigma(x,\xi)$.
}

More important point is that symbol of an operator should be more simple to verify its Fredholm properties. For two above examples a Fredholm property of an operator symbol is equivalent to its invertibility.

\section{Stratification of manifolds and operators}

\subsection{Sub-manifolds}

The above definition of an operator on a manifold supposes that all neighborhoods $\{U_j\}$ have the same type. But if a manifold has a smooth boundary even then there are two types of neighborhoods related to a place of neighborhood, namely inner neighborhoods and boundary ones. For inner neighborhood $U$ such that $\overline{U}\subset\mathring{M}$ we have the diffeomorphism $\omega: U\rightarrow D$, where $D\in\mathbb R^m$ is an open set.                                                                                                                                                                                                                                                                                                                                                                                                                                                                                                                                                                                                                                                                                                                                                                                                                                                                                                                                                                                                                                                                                                                                                                                                                                                                                                                                                                                                                                                                                                                                                                                                                                                                                                                                                                                                                                                                                                                                                                                                                                                                                                                                                                                                                                                                                                                                                                                                                                                                                                                                                                  For a boundary neighborhood such that $U\cap\partial M\neq\emptyset$ we have another diffeomorphism  $\omega_1: U\rightarrow D\cap\mathbb R^m_+$, where $\mathbb R^m_+=\{x\in\mathbb R^m: x=(x_1,\cdots,x_m), x_m>0\}$. May be this boundary $\partial M$ has some singularities like conical points and wedges. The conical point at the boundary is a such point, for which its neighborhood is diffeomorphic to the cone
$
C^a_+=\{x\in{\mathbb R}^m:~x_m>a|x'|,~x'=(x_1,...,x_{m-1}),~a>0\},
$
the wedge point  of codimension $k, 1\leq k\leq m-1,$ is a such point for which its neighborhood is diffeomorphic to the set   $\{x\in\mathbb R^m: x=(x',x''), x''\in\mathbb R^{m-k}, x'=(x_1,\cdots,x_{m-k-1}), x_{m-k-1}>a|x'''|, x'''=(x_1,\cdots,x_{m-k-2}), a>0\}$. So if the manifold $M$ has such singularities we suppose that we can extract certain k-dimensional sub-manifolds, namely $(m-1)$-dimensional boundary $\partial M$,  and $k$-dimensional wedges $M_k, k=0,\cdots,m-2$; $M_0$ is a collection of conical points.

\subsection{Enveloping operators}

If the family $\{A_x\}_{x\in M}$ is continuous in operator topology, then according to Simonenko's theory there is an enveloping operator, i.e. such an operator $A$ for which every operator $A_x$ is the local representative for the operator $A$ in the point $x\in M$.

{\bf Example 5.} {\it
If $\{A_x\}_{x\in M}$ consists of Calderon--Zygmund operators in $\mathbb R^m$ \cite{MPr} with symbols $\sigma_x(\xi)$ parametrized by points $x\in M$ and this family smoothly depends on $x\in M$ then Calderon--Zygmubd operator with variable kernel and symbol $\sigma(x,\xi)$ will be an enveloping operator for this family.
}

{\bf Example 6.} {\it
If $\{A_x\}_{x\in M}$ consists of null operators then an enveloping operator is a compact operator \cite{S}.
}

{\bf Theorem 1.} {\it
Operator $A$ has a Fredholm property iff its all local representatives $\{A_x\}_{x\in M}$ have the same property.
}

This property was proved in \cite{S}. But we will give the proof (see Lemma 2) including some new constructions because it will be used below for a decomposition of the operator.

\subsection{Hierarchy of operators}

We will remind here the following definition and Fredholm criteria for operators \cite{K}.

{\bf Definition 8.} {\it
Let $B_1, B_2$ be Banach spaces, and $A: B_1\rightarrow B_2$ be a linear bounded operator. The operator $R: B_2\rightarrow B_1$ is called a regularizer for the operator $A$ if the following properties
\[
RA=I_1+T_1,~~~AR=I_2+T_2
\]
hold, where $I_k: B_k\rightarrow B_k$ is identity operator, $T_k; B_k\rightarrow B_k$ is a compact operator, $k=1,2$.
}

{\bf Proposition 2.} {\it The operator $A: B_1\rightarrow B_2$ has a Fredholm property iff there exists a linear bounded regularizer $R: B_2\rightarrow B_1$.
}

{\bf Lemma 1.} {\it
Let $f$ be a smooth function on the manifold $M, U\subset M$ be an open set, and $supp~f\subset U$. Then the operator $f\cdot A-A\cdot f$ is a compact operator.
}

\begin{proof}
Let $g$ be a smooth function on $M, supp~g\subset V\subset M$, moreover $\overline U\subset V, g(x)\equiv 1$ for $x\in supp~f$. Then we have
\[
f\cdot A=f\cdot A\cdot g+f\cdot A\cdot(1-g)=f\cdot A\cdot g+T_1,
\]
\[
A\cdot f=g\cdot A\cdot f+(1-g)\cdot A\cdot f=g\cdot A\cdot f+T_2,
\]
where $T_1,T_2$ are compact operators.  Let us denote $g\cdot A\cdot g\equiv h$ and write
\[
f\cdot A\cdot g=f\cdot g\cdot A\cdot g=f\cdot h,~~~g\cdot A\cdot f=g\cdot A\cdot g\cdot f=h\cdot f,
\]
and we obtain the required property.
\end{proof}

{\bf Definition 9.} {\it The operator $A$ is called an elliptic operator if its operator symbol $\{A_x\}_{x\in M}$ consists of Fredholm operators.
}

Now we will show that each elliptic operator really has a Fredholm property. Our proof in general follows the book \cite{S}, but our constructions are more stratified and we need such constructions below.

{\bf Lemma 2.} {\it
Let $A$ be an elliptic operator. Then the operator $A$ has a Fredholm property.
}

\begin{proof}
To obtain the proof we will construct the regularizer for the operator $A$. For this purpose we choose two coverings like proposition 1 and write the operator $A$ in the form
\begin{equation}\label{1}
A=\sum\limits_{j=1}^{n}f_j\cdot A\cdot g_j +T,
\end{equation}
where $T$ is a compact operator. Without loss of generality we can assume that there are $n$ points $x_k\in U_k\subset V_k, k=1,2,...n$. Moreover, we can construct such coverings by balls in the following way. Let $\varepsilon>0$ be enough small number. First, for every point $x\in M_0$ we take two balls $U_x,V_x$ with the center at $x$ of radius $\varepsilon$ and $2\varepsilon$ and construct two open coverings for $M_0$ namely $\frak U_0=\cup_{x\in M_0}U_x$ and $\frak V_0=\cup_{x\in M_0}V_x$ . Second, we consider the set $L_1=\overline M\setminus\frak V_0$ and construct two coverings $\frak U_1=\cup_{x\in L_1\cap M_1}U_x$ and $\frak V_1=\cup_{x\in L_1\cap M_1}V_x$. Further, we introduce the set $L_2=\overline M\setminus(\frak V_0\cup\frak V_1)$ and two coverings  $\frak U_2=\cup_{x\in L_2\cap M_2}U_x$ and $\frak V_2=\cup_{x\in L_2\cap M_2}V_x$. Continuing these actions we will come to the set $L_{m-1}=\overline M\setminus(\cup_{k=0}^{m-2}\frak U_k)$ which consists of smoothness points of $\partial M$ and inner points of $M$, construct two covering $\frak U_{m-1}=\cup_{x\in L_{m-1}\cap\partial M}U_x$ and $\frak V_{m-1}=\cup_{x\in L_{m-1}\cap\partial M}V_x$.
Finally, the set $L_m=\overline M\setminus(\cup_{k=0}^{m-1}\frak U_k)$ consists of inner points of the manifold $M$ only. We finish this process by choosing the covering $\frak U_m$ for the latter set $L_m$. So, the covering $\cup_{k=0}^m\frak U_k$ will be a covering for the whole manifold $M$. According to compactness property we can take into account that this covering is finite, and centers of balls which cover $M_k$ are placed at $M_k$.

Now we will rewrite the formula \eqref{1} in the following way
\begin{equation}\label{4}
A=\sum\limits_{k=0}^m\left(\sum\limits_{j=1}^{n_k}f_{jk}\cdot A\cdot g_{jk}\right) +T,
\end{equation}
where coverings and partitions of unity $\{f_{jk}\}$ and $\{g_{jk}\}$ are chosen as mentioned above. In other words the operator
\[
\sum\limits_{j=1}^{n_k}f_{jk}\cdot A\cdot g_{jk}
\]
is related to some neighborhood of the sub-manifold $M_k$; this neighborhood is generated by covering of the sub-manifold $M_k$ by balls with centers at points $x_{jk}\in M_k$. Since $A_{x_{jk}}$ is a local representative for the operator $A$ at point $x_{jk}$ we can rewrite the formula \eqref{4} as follows
\begin{equation}\label{5}
A=\sum\limits_{k=0}^m\left(\sum\limits_{j=1}^{n_k}f_{jk}\cdot A_{x_{jk}}\cdot g_{jk}\right) +T.
\end{equation}

Let us denote $S_{\omega_j^{-1}}\hat g_j\equiv\tilde g_j, \hat f_jS_{\omega_j}\equiv\tilde f_j$.
Further, we can assert that the operator
\[
R=\sum\limits_{k=0}^m\left(\sum\limits_{j=1}^{n_k}g_{jk}\cdot A^{-1}_{x_{jk}}f_{jk}\right),
\]
will be the regularizer for the operator $A'$; here $ A^{-1}_{x_{jk}}$ is a regularizer for the operator $A_{x_{jk}}$

Indeed,
\[
RA=\left(\sum\limits_{k=0}^m\left(\sum\limits_{j=1}^{n_k}g_{jk}A_{x_{jk}}^{-1}f_{jk}\right)\right)\cdot A=
\]
\[
\sum\limits_{k=0}^m\sum_{j=1}^{n_k}g_{jk}\cdot A^{-1}_{x_{jk}}\cdot (A-A_{x_{jk}}+A_{x_{jk}})\cdot f_{jk}+T_1=
\]
\[
\sum\limits_{k=0}^m\sum_{j=1}^{n_k}g_{jk}\cdot A^{-1}_{x_{jk}}\cdot (A-A_{x_{jk}})\cdot f_{jk}+\sum\limits_{k=0}^m\sum_{j=1}^{n_k}f_{jk}+T_1=I_1+T_1+\Theta_1,
\]
\[
\Theta_1=\sum\limits_{k=0}^m\sum_{j=1}^{n_k}g_{jk}\cdot A^{-1}_{x_{jk}}\cdot (A-A_{x_{jk}})\cdot f_{jk},
\]
because $f_{jk}\cdot A_{x_{jk}}= A_{x_{jk}}\cdot f_{jk}+\text{compact summand}$, and $f_{jk}\cdot g_{jk}=f_{jk}$, and
\[
\sum\limits_{k=0}^m\sum_{j=1}^{n_k}f_{jk}\equiv 1
\]
as the partition of unity. The same property
\[
AR=I_2+T_2+\Theta_2,
\]
\[
\Theta_2=\sum\limits_{k=0}^m\sum_{j=1}^{n_k}g_{jk}\cdot (A-A_{x_{jk}})\cdot A^{-1}_{x_{jk}}\cdot f_{jk},
\]
is verified analogously.
\end{proof}

\section{Piece-wise continuous operator families}

Given operator $A$ with the symbol $\{A_x\}_{x\in\overline M}$ generates a few operators in dependence on a quantity of singular manifolds. We consider this situation in the following way. We will assume additionally some smoothness properties for the symbol $\{A_x\}_{x\in\overline M}$.

{\bf Theorem 2.} {\it
If the symbol $\{A_x\}_{x\in\overline M}$ is a piece-wise continuous operator function then there are $m+1$ operators $A^{(k)}, k=0,1,\cdots,m$ such that the operator $A$ and the operator
\begin{equation}\label{2}
A'=\sum\limits_{k=0}^mA^{(k)}+T
\end{equation}
have the same symbols, where the operator $A^{(k)}$ is an enveloping operator for the family $\{A_x\}_{x\in\overline M_k}, T$ is a compact operator.
}

\begin{proof}
We will use the constructions from proof of Lemma 2, namely the formula \eqref{5}. We will extract the operator
\[
\sum\limits_{j=1}^{n_k}f_{jk}\cdot A_{x_{jk}}\cdot g_{jk}
\]
which ``serves'' the sub-manifold $M_k$ and consider it in details. This operator is related to neighborhoods $\{U_{jk}\}$ and the partition of unity $\{f_{jk}\}$. Really $U_{jk}$ is the ball with the center at $x_{jk}\in M_k$ of radius $\varepsilon>0$, therefore $f_{jk}, g_{jk}, n_k$ depend on $\varepsilon$.

According to Simonenko's ideas \cite{S} we will construct the component $A^{(k)}$ in the following way. Let $\{\varepsilon_n\}_{n=1}^{\infty}$  be a sequence such that $\varepsilon_n>0, \forall n\in\mathbb N, \lim\limits_{n\to\infty}\varepsilon_n=0$. Given $\varepsilon_n$ we choose coverings $\{U_{jk}\}_{j=1}^{n_k}$ and $\{V_{jk}\}_{j=1}^{n_k}$ as above with partition of unity $\{f_{jk}\}$ and corresponding functions $\{g_{jk}\}$ such that
\[
|||f_{jk}\cdot(A_x-A_{x_{jk}})\cdot g_{jk}|||<\varepsilon_n,~~~\forall x\in V_{jk};
\]
we remaind that $U_{jk}, V_{jk}$   are balls with centers at $x_{jk}\in\overline{M_k}$ of radius $\varepsilon$ and $2\varepsilon$. This requirement is possible according to continuity of family $\{A_x\}$ on the sub-manifold $\overline{M_k}$. Now we will introduce such constructed operator
\[
A_n=\sum\limits_{j=1}^{n_k}f_{jk}\cdot A_{x_{jk}}\cdot g_{jk}
\]
and will show that the sequence $\{A_n\}$ is Cauchy sequence with respect to norm $|||\cdot|||$. We have
\[
A_l=\sum\limits_{i=1}^{l_k}F_{ik}\cdot A_{y_{ik}}\cdot G_{ik},
\]
where the operator $A_l$ is constructed for given $\varepsilon_l$ with corresponding coverings $\{u_{ik}\}_{i=1}^{l_k}$ and $\{v_{ik}\}_{j=1}^{i_k}$ with partition of unity $\{F_{ik}\}$ and corresponding functions $\{G_{ik}\}$ so that
\[
|||F_{ik}\cdot(A_x-A_{y_{ik}})\cdot G_{ik}|||<\varepsilon_l,~~~\forall x\in v_{ik};
\]
here $u_{ik}, v_{ik}$   are balls with centers at $y_{ik}\in\overline{M_k}$ of radius $\tau$ and $2\tau$.

We can write
\[
A_n=\sum\limits_{j=1}^{n_k}f_{jk}\cdot A_{x_{jk}}\cdot g_{jk}=\sum\limits_{i=1}^{l_k}F_{ik}\cdot\sum\limits_{j=1}^{n_k}f_{jk}\cdot A_{x_{jk}}\cdot g_{jk}=
\]
\[
\sum\limits_{i=1}^{l_k}\sum\limits_{j=1}^{n_k}F_{ik}\cdot f_{jk}\cdot A_{x_{jk}}\cdot g_{jk}=\sum\limits_{i=1}^{l_k}\sum\limits_{j=1}^{n_k}F_{ik}\cdot f_{jk}\cdot A_{x_{jk}}\cdot g_{jk}\cdot G_k+T_1,
\]
and the same we can write for $A_l$
\[
A_l=\sum\limits_{i=1}^{l_k}F_{ik}\cdot A_{y_{ik}}\cdot G_{ik}=\sum\limits_{j=1}^{n_k}f_{jk}\cdot \sum\limits_{i=1}^{l_k}F_{ik}\cdot A_{y_{ik}}\cdot G_{ik}=
\]
\[
\sum\limits_{j=1}^{n_k}\sum\limits_{i=1}^{l_k}f_{jk}\cdot F_{ik}\cdot A_{y_{ik}}\cdot G_{ik}=\sum\limits_{j=1}^{n_k}\sum\limits_{i=1}^{l_k}f_{jk}\cdot F_{ik}\cdot A_{y_{ik}}\cdot G_{ik}\cdot g_{jk}+T_2.
\]

Let us consider the difference
\begin{equation}\label{6}
|||A_n-A_l|||=|||\sum\limits_{j=1}^{n_k}\sum\limits_{i=1}^{l_k}f_{jk}\cdot F_{ik}\cdot (A_{x_{jk}}-A_{y_{ik}})\cdot G_{ik}\cdot g_{jk}|||.
\end{equation}

Obviously, summands with non-vanishing supplements to the formula \eqref{6} are those for which $U_{jk}\cap u_{ik}\neq\emptyset$. A number of such neighborhoods are finite always for arbitrary finite coverings, hence we obtain
\[
|||A_n-A_l|||\leq\sum\limits_{j=1}^{n_k}\sum\limits_{i=1}^{l_k}|||f_{jk}\cdot F_{ik}\cdot (A_{x_{jk}}-A_{y_{ik}})\cdot G_{ik}\cdot g_{jk}|||\leq
\]
\[
\sum\limits_{x\in U_{jk}\cap u_{ik}\neq\emptyset}|||f_{jk}\cdot F_{ik}\cdot (A_{x_{jk}}-A_x)\cdot G_{ik}\cdot g_{jk}|||+
\]
\[
\sum\limits_{x\in U_{jk}\cap u_{ik}\neq\emptyset}|||f_{jk}\cdot F_{ik}\cdot (A_x-A_{y_{ik}})\cdot G_{ik}\cdot g_{jk}|||\leq 2K\max[\varepsilon_n,\varepsilon_l],
\]
where $K$ is a universal constant.

Thus, we have proved that the sequence $\{A_n\}$ is a Cauchy sequence hence there exists $\lim\limits_{n\to\infty}A_n=A^{(k)}$.
\end{proof}

{\bf Corollary 1.} {\it
The operator $A$ has a Fredholm property iff all operators $A^{(k)}, k=0,1,\cdots,m$ have the same property.
}

{\bf Remark 2.} {\it
The constructed operator $A'$ generally speaking does not coincide with the initial  operator $A$ because they act in different spaces. But for some cases they may be the same.
}

\section{Conclusion}

This paper is a general concept of my vision to the theory of pseudo-differential equations and boundary value problems on manifolds with a non-smooth boundary. The second part will be devoted to applying these abstract results to index theory for such operator families and then to concrete classes of pseudo-differential equations.

\end{document}